# On the Łojasiewicz exponent at infinity for polynomial functions

Laurenţiu Păunescu , Alexandru Zaharia

## 1 Introduction

**1.1** For $n, q \in \mathbf{N} \setminus \{0\}$ we consider the polynomial functions

$$f = f_{n,q} : \mathbf{C}^3 \longrightarrow \mathbf{C} \quad , \quad f(x,y,z) = f_{n,q}(x,y,z) := x - 3x^{2n+1}y^{2q} + 2x^{3n+1}y^{3q} + yz \quad .$$

We will study some properties of these polynomials , related to their behaviour at infinity , and we will prove that some results , obtained in [14] and [3] , [4] , for the case of polynomials in two variables , are not true in the case of polynomials in $m \geq 3$ variables . Also , our polynomials $f_{n,q}$ show that several classes of polynomials , with "good" behaviour at infinity , considered in [6] , [7] , [12] , [10] , are distinct .

The first remark on our polynomials is :

**1.2 Remark .** After a suitable polynomial change of coordinates in $\mathbf{C}^3$ , one can write $f(X,y,Z) = X$ . Namely , taking $Z := z - 3x^{2n+1}y^{2q-1} + 2x^{3n+1}y^{3q-1}$ , we get : $f(x,y,Z) = x + yZ$ . Next , we put $X := x + yZ$ and we obtain $f(X,y,Z) = X$ . Thus , there exists a polynomial automorphism $P = (P_1, P_2, P_3) : \mathbf{C}^3 \longrightarrow \mathbf{C}^3$ such that $f = P_1$ .

**1.3** For a polynomial $g : \mathbf{C}^m \longrightarrow \mathbf{C}$ , we consider grad $g(x) := \left( \overline{\frac{\partial f}{\partial x_1}(x)}, \ldots, \overline{\frac{\partial f}{\partial x_m}(x)} \right)$ . If $g$ has non–isolated singularities , the *Łojasiewicz number at infinity* , $L_\infty(g)$ , is defined by $L_\infty(g) := -\infty$ . When $g$ has only isolated singularities , the Łojasiewicz number at infinity is the supremum of the set

$$\{ \ \nu \in \mathbf{R} \ | \ \exists A > 0 \ , \ \exists B > 0 \ , \ \forall x \in \mathbf{C}^m \quad , \quad \text{if} \ \ \|x\| \geq B \quad , \quad \text{then} \ \ A\|x\|^\nu \leq \|\text{grad } g(x)\| \ \} \quad .$$

Equivalent definition is ( see for instance [14] or [5] , proof of Proposition 1 ) :

$$L_\infty(g) := \lim_{r \to \infty} \frac{\log \varphi(r)}{\log r} \quad , \quad \text{where} \ \ \varphi(r) := \inf_{\|x\|=r} \|\text{grad } g(x)\| \quad .$$

The following result is a reformulation of Theorem 10.2 from [4] :

**Theorem .** *Let* $g : \mathbf{C}^2 \longrightarrow \mathbf{C}$ *be a polynomial function . Then there exists a polynomial automorphism* $P = (P_1, P_2) : \mathbf{C}^2 \longrightarrow \mathbf{C}^2$ *such that* $g = P_1$ *if and only if* $g$ *has no critical values and* $L_\infty(g) > -1$ .

**1.4** In the next Section we will prove the following

**Proposition .** $L_\infty(f_{n,q}) = -\dfrac{n}{q}$ .

In particular , if $n \geq q$ , then $L_\infty(f_{n,q}) \leq -1$ . Using Remark 1.2 , our Proposition shows that Theorem 1.3 can not be extended to the case of a polynomial function $g : \mathbf{C}^m \longrightarrow \mathbf{C}$ , when $m \geq 3$ .

**1.5** It is proved in [14] and [4] that a polynomial function $g : \mathbf{C}^2 \longrightarrow \mathbf{C}$ has $L_\infty(g) \neq -1$ . Proposition 1.4 shows that this is no longer true for polynomial functions $\mathbf{C}^m \longrightarrow \mathbf{C}$ , when $m \geq 3$ .

**1.6** If $g : \mathbf{C}^m \longrightarrow \mathbf{C}$ is a polynomial function , we call $t_0 \in \mathbf{C}$ a *typical* value of $g$ if there exists $U \subseteq \mathbf{C}$ , an open neigbourhood of $t_0$ , such that the restricition $g : g^{-1}(U) \longrightarrow U$ is a $C^\infty$ trivial fibration . If $t_0$ is not a typical value of $g$ , then $t_0$ is called an *atypical* value of $g$ . In general , the *bifurcation set* , $B_g$ , of atypical values of $g$ , contains , besides the set $\Sigma_g$ of critical values og $g$ , some extra values , the so–called "critical values coming from infinity" . For example , if $g(x,y) = x^2y + x$ , then $\Sigma_g = \emptyset$ and $B_g = \{0\}$ .

Several classes of polynomials without "critical values coming from infinity" are considered in literature , see for instance [1] , [6] , [7] , [11] , [9] . We recall now three of them . In the next Section we will use the polynomials $f_{n,q}$ to show that these classes are distinct .

For a polynomial $g : \mathbf{C}^m \longrightarrow \mathbf{C}$ , we denote :

$$M(g) := \{\ x \in \mathbf{C}^m \ \mid \exists \lambda \in \mathbf{C} \ \text{ such that } \ \text{grad}\ g(x) = \lambda \cdot x\ \} \quad .$$

Geometrically , a point $x \in M(g)$ if and only if either $x$ is a critical point of $g$ , or $x$ is not a critical point of $g$ , but the hypersurface $g^{-1}(g(x))$ does not intersect transversally , at $x$ , the sphere $\{\ z \in \mathbf{C}^m \mid \|z\| = \|x\|\ \}$ .

A polynomial $g : \mathbf{C}^m \longrightarrow \mathbf{C}$ is called *M–tame* if for any sequence $\{z^k\} \subseteq M(g)$ such that $\lim\limits_{k\to\infty} \|z^k\| = \infty$ , we have $\lim\limits_{k\to\infty} |g(z^k)| = \infty$ . See [10] for properties of M–tame polynomials .

A polynomial $g : \mathbf{C}^m \longrightarrow \mathbf{C}$ is called *quasitame* if for any sequence $\{z^k\} \subseteq \mathbf{C}^m$ such that $\lim\limits_{k\to\infty} \|z^k\| = \infty$ and $\lim\limits_{k\to\infty} \text{grad}\ g(z^k) = 0$ , we have $\lim\limits_{k\to\infty} |g(z^k) - \langle z^k, \text{grad}\ g(z^k)\rangle| = \infty$ . Here , $\langle\cdot,\cdot\rangle$ denotes the Hermitian product on $\mathbf{C}^m$ . See [6] , [7] for properties of quasitame polynomials .

Follwing [11] , [12] , we will say that a polynomial $g : \mathbf{C}^m \longrightarrow \mathbf{C}$ satisfies Malgrange's condition for $t_0 \in \mathbf{C}$ if , for $\|x\|$ large enough and for $g(x)$ close to $t_0$ , there exists $\delta > 0$ such that $\|x\| \cdot \|\text{grad}\ g(x)\| \geq \delta$ . Equivalent formulation is : there exists no sequence $\{z^k\} \subseteq \mathbf{C}^m$ such that $\lim\limits_{k\to\infty} \|z^k\| = \infty$ , $\lim\limits_{k\to\infty} g(z^k) = t_0$ and $\lim\limits_{k\to\infty} \|z^k\| \cdot \|\text{grad}\ g(z^k)\| = 0$ .

The next result seems to be well–known , see [10] , [8] , [13] . Its proof can be easily obtained , by contradiction .

**1.7 Proposition .** *For* $m \geq 2$ *, let* $g : \mathbf{C}^m \longrightarrow \mathbf{C}$ *be a polynomial function .*

($a$) *If* $g$ *is quasitame , then* $g$ *satisfies Malgrange's condition for any* $t_0 \in \mathbf{C}$ *.*

($b$) *If* $g$ *satisfies Malgrange's condition for any* $t_0 \in \mathbf{C}$ *, then* $g$ *is M–tame .*

We will show that these implications can not be reversed , if $m \geq 3$ ( see also [2] ) . More precisely , we have :

**1.8 Proposition .** *(a) For any $n, q \in \mathbf{N} \setminus \{0\}$ , the polynomial $f_{n,q}$ is M–tame , but not quasitame .*

*(b) The polynomial $f_{n,q}$ satisfies Malgrange's condition for any $t_0 \in \mathbf{C}$ , if and only if $n \leq q$ .*

Thus , if $f = f_{n,q}$ for some $n > q$ , then , by [12] , the family $\overline{f}$ of projective closures of fibres of $f$ has nontrivial vanishing cycles , despite Remark 1.2 . Also , such an $f$ is not $t$–regular at infinity , in the sense of [13] , since by [12] , the $t$–regularity at infinity is equivalent to Malgrange's condition for any $t_0 \in \mathbf{C}$ .

# 2 Proofs

**2.1** Let $g : \mathbf{C}^m \longrightarrow \mathbf{C}$ be a polynomial function with only isolated singularities . For an analytic curve $p : (0, \varepsilon) \longrightarrow \mathbf{C}^m$ such that $\lim\limits_{t \to \infty} \|p(t)\| = \infty$ , we consider the expansions in Laurent series :

$$p(t) = at^{\alpha} + a_1 t^{\alpha+1} + \dots \quad , \quad \text{with} \;\; \alpha < 0 \;\; \text{and} \;\; a \neq 0 \tag{1}$$

$$\text{grad } g(p(t)) = bt^{\beta} + b_1 t^{\beta+1} + \dots \quad , \quad \text{with} \;\; \beta \neq 0 \;\; \text{and} \;\; b \neq 0 \tag{2}$$

and we denote $L(g;p) := \dfrac{ord(\text{grad } g(p(t))}{ord(p(t))} = \dfrac{\beta}{\alpha}$ . Here , $ord$ denotes the order of series . It follows , for example from ( the proof of ) Proposition 1 in [5] , that

$$L_{\infty}(g) = \inf \left\{ L(g;p) \;\middle|\; \begin{array}{l} p : (0, \varepsilon) \longrightarrow \mathbf{C}^m \text{ is an analytic curve} \\ \text{such that (1) and (2) are fulfilled} \end{array} \right\} \tag{3}$$

**2.2 Proof of Proposition 1.4 .** Consider the curve $\Psi : (0,1) \longrightarrow \mathbf{C}^m$ defined by : $\Psi(t) := (t^{-q}, t^n, 0)$ . Then $L(f_{n,q}; \Psi) = -\dfrac{n}{q}$ , hence $L_{\infty}(f_{n,q}) \leq -\dfrac{n}{q} < 0$ .

Let now consider an arbitrary analytic curve $p : (0, \varepsilon) \longrightarrow \mathbf{C}^m$ , $p(t) = (x(t), y(t), z(t))$ , such that $\lim\limits_{t \to 0} \|p(t)\| = \infty$ .

If $\lim\limits_{t \to 0} \|\text{grad } f_{n,q}(p(t))\| \neq 0$ , then $ord(\text{grad } f_{n,q}(p(t))) \leq 0$ . Hence , $L(f_{n,q}; p) \geq 0 > L(f_{n,q}; \Psi)$ .

Suppose now that $\lim\limits_{t \to 0} \|\text{grad } f_{n,q}(p(t))\| = 0$ . Then $\lim\limits_{t \to 0} y(t) = 0$ and

$$\lim_{t \to 0} (x(t))^n \cdot (y(t))^q \;\; \text{is a root of the equation} \;\; 1 - (6n+3)T^2 + (6n+2)T^3 = 0 \quad . \tag{4}$$

Hence

$$y(t) \not\equiv 0 \quad , \quad \lim_{t \to 0} y(t) = 0 \;\; \text{and} \;\; \lim_{t \to 0} \|x(t)\| = \infty \quad . \tag{5}$$

Therefore , we have :

$$ord(y(t)) \geq ord(\text{grad } f_{n,q}(p(t))) > 0 \;\; \text{and} \;\; ord(p(t)) \leq ord(x(t)) < 0 \quad . \tag{6}$$

It follows , using (4) and (6) , that

$$L(f_{n,q};p) = \frac{ord(\text{grad } f_{n,q}(p(t)))}{ord(p(t))} \geq \frac{ord(y(t))}{ord(p(t))} \geq \frac{ord(y(t))}{ord(x(t))} = -\frac{n}{q} \quad .$$

Proposition 1.4 is proved .

◊

**2.3 Proof of Proposition 1.8 .** Part $(b)$ follows from ( the proof of ) Proposition 1.4 .

If $n, q \in \mathbf{N} \setminus \{0\}$ are fixed , then the curve $\Psi(t) := (t^{-q}, t^n, 0)$ can be used to show that $f_{n,q}$ is not quasitame .

Suppose now that $f = f_{n,q}$ is not M–tame . Using Curve Selection Lemma at infinity , see [10] , one can find an analytic curve $p : (0, \varepsilon) \longrightarrow M(f)$ , $p(t) = (x(t), y(t), z(t))$ , such that $\lim_{t\to 0} f(p(t)) \in \mathbf{C}$ . This implies that $\lim_{t\to 0} \text{grad } f(p(t)) = 0$ , hence relations (4) and (5) hold . The condition $p(t) \in M(f)$ means that

$$\left(\overline{\frac{\partial f}{\partial x}(p(t))} \ , \ \overline{\frac{\partial f}{\partial y}(p(t))} \ , \ \overline{\frac{\partial f}{\partial z}(p(t))}\right) = \lambda(t) \cdot (x(t), y(t), z(t)) \tag{7}$$

for some suitable analytic curve $\lambda : (0, \varepsilon) \longrightarrow \mathbf{C}$ . It follows that none of the components of $p(t)$ or of grad $f(p(t))$ is identically zero .

If $A := ord(x(t))$ , $B := ord(y(t))$ and $C := ord\left(\frac{\partial f}{\partial x}(p(t))\right)$ , then $A < 0$ , $B > 0$ , $C > 0$ , and relations (4) and (7) give us

$$ord(\lambda(t)) = C - A \quad , \quad nA + qB = 0 \ \text{ and } \ ord(z(t)) = B + A - C \quad . \tag{8}$$

Since $f = y\frac{\partial f}{\partial y} + x\left(1 + (6q-3)x^{2n}y^{2q} - (6q-2)x^{3n}y^{3q}\right)$ , it is easy to see that

$$\lim_{t\to 0}(x(t))^n \cdot (y(t))^q \ \text{ is a root of the equation } \ 1 + (6q-3)T^2 - (6q-2)T^3 = 0 \quad .$$

Thus , using (4) , it follows that $\lim_{t\to 0}(x(t))^n \cdot (y(t))^q = 1$ . Hence we can assume that

$$x(t) = t^A \quad , \quad y(t) = t^B + t^{B+D} \cdot \rho(t) \quad , \quad \text{with } \ D > 0 \ \text{ and } \ \rho(0) \neq 0$$

( since $(x(t))^n \cdot (y(t))^q \equiv 1$ implies that $\frac{\partial f}{\partial x}(p(t)) \equiv 0$ ) . By a direct computation , we find that

$$\lambda(t) = \frac{1}{x(t)} \cdot \overline{\frac{\partial f}{\partial x}(p(t))} = 6nq \cdot \overline{\rho(0)} \cdot t^{D-A} + \text{higher terms} \quad ,$$

hence , by (8) , $D = C$ . Next , the second component in (7) gives us

$$z(t) = 6qx(t)^{2n+1}y(t)^{2q-1}\left(1 - x(t)^n y(t)^q\right) + \overline{\lambda(t)y(t)} =$$

$$= -6q^2 \cdot \rho(0) \cdot t^{A-B+D} + \text{higher terms} \quad .$$

Hence , by (8) , we have $B = C = D$ . Finally , the third component in (7) gives us

$$t^B + t^{B+D} \cdot \overline{\rho(t)} = 6nq \cdot \overline{\rho(0)} \cdot t^{D-A} \cdot \left(-6q^2 \cdot \rho(0) \cdot t^{A-B+D}\right) + \text{higher terms} \quad .$$

Comparing the ledings coefficients , we obtain

$$1 = -36nq^3 \cdot |\rho(0)|^2 \quad ,$$

which is impossible . Thus , Proposition 1.8 is proved .

$\diamondsuit$

**2.4 Remark .** (*i*) If $n \leq q$ , it is possible to prove that $f_{n,q}$ is M–tame just looking at the order of various Laurent expansions .

(*ii*) It is not difficult to see that for the polynomials $f_{n,q}$ , the Newton nondegeneracy condition fails on a face of dimension 1 . Using this , one can construct other similar examples .

**Aknowledgements .** *The second author would like to thank The Fields Institute and the University of Toronto , where this note was written , for hospitality and support .*

Laurenţiu Păunescu
University of Sydney
School of Mathematics and Statistics
Sydney
NSW 2006
Australia

Alexandru Zaharia
Institute of Mathematics of
the Romanian Academy
*Mailing address* ( *until July* 1997 ) :
The Fileds Institute
222 College Street
Toronto
Ontario , M5T 3J1
Canada